\documentstyle[12pt, twoside]{article}
\input amssym.def
\input amssym.tex

\newenvironment{demo}[1]%
{\vskip-\lastskip\medskip
  \noindent
  {\em #1.}\enspace
  }%
{\qed\par\medskip
  }

\newcommand{\qed}{
  \strut\hfill
  \mbox{$\Box$}
  }

\newtheorem{theorem}{Theorem}

\newtheorem{remark}{Remark}

\newtheorem{proposition}{Proposition}

\newcommand{\bdl}{\underline{{\C}^r}}
\newcommand{\C}{\Bbb C}
\newcommand{\Cg}{{\cal C}_{\lambda} }
\newcommand{\Cbar}{\overline{\cal C}_{\lambda} }
\newcommand{\Ct}{\widetilde{\cal C}_{\lambda}}
\newcommand{\F}{{\cal F} (m, n -s, V)}
\newcommand{\Fm}{{\cal F} (m, n -m, V)}
\newcommand{\Gr}{Gr(m, V)}
\newcommand{\Grdual}{Gr(n -s, V)}
\newcommand{\Grdualm}{Gr(n -m, V)}

\newcommand{\Hom}{\mbox{Hom}}
\newcommand{\NCQ}{{\cal{NCQ}\/} }
\newcommand{\NCQT}{\widetilde{\cal{NCQ}} }
\newcommand{\Nt}{\widetilde{\cal N}}
\newcommand{\Ntt}{\widetilde{\cal N}_1}
\newcommand{\Ntd}{\widetilde{\cal N}_2}
\newcommand{\Obar}{\overline{\cal O}_{\lambda} }
\newcommand{\Ot}{\widetilde{\cal O}_{\lambda}}
\newcommand{\Ott}{\widetilde{\cal O}_{\lambda}^1}
\newcommand{\Otd}{\widetilde{\cal O}_{\lambda}^2}
\newcommand{\Taut}{ \mbox{Taut} }

\begin{document}

\title{Resolution of singularities of null cones
}
\author{Weiqiang Wang\\ \small
Department of Mathematics\\\small North Carolina State University \\
\small Raleigh, NC 27695\\ \small wqwang@math.ncsu.edu}

%words{Null cone; resolution of singularities; Lie groups}

\date{}

%%% ----------------------------------------------------------------------

%%% ----------------------------------------------------------------------
\maketitle
\begin{abstract}
 We give canonical resolutions of singularities of several
 cone varieties arising from invariant theory.
 We establish a connection between our resolutions and
 resolutions of singularities
 of closure of conjugacy classes in classical Lie algebras.
\end{abstract}
%%% ----------------------------------------------------------------------

\section*{Introduction}
  The purpose of this paper is to present canonical resolutions
of singularities of certain cone varieties arising naturally
in invariant theory. Examples of such cone varieties are provided by the
so-called null cones which are the zero locus of the $O_n$ (or $Sp_n$)
invariant polynomials on $M_{n,m}$,
where $M_{n,m}$ denotes the space of $n \times m$ complex matrices
(cf. e.g. \cite{H}).
Another example is the null cone for general linear groups.
Our construction resembles
the celebrated Springer resolution of the nilpotent variety
in a complex semisimple Lie algebra
which plays a vital role
in representation theory and geometry (cf. \cite{CG}).

It turns out that the cone varieties studied here include
as special cases the
variety $Z$ constructed in \cite{KP1, KP2} associated to
two-column partitions. The Kraft-Procesi variety $Z$
is a complete intersection whose
quotient under certain group coincide with the closure
of a nilpotent conjugacy class in classical Lie algebras.
It is a very interesting problem to construct a canonical
resolution of singularities of $Z$ in general.
We establish relations between our resolutions of singularities
and those for closure of conjugacy classes associated
to two-column partitions. In the general linear
Lie algebra case, we find more than one canonical
resolution of singularities of null cones and of
closure of conjugacy classes.

We mention in passing that
the present work grows out of an attempt to generalize
the geometric construction of $(gl_n, gl_m)$-duality \cite{W}
in the spirit of \cite{CG} to other Howe duality \cite{H}.

The paper is divided into two sections. In Sect.~\ref{sect_classic}
we present resolutions of singularities of null cones in the
$O_n$ and $Sp_n$ cases, and establish connections
with resolutions of singularities of closure of conjugacy
classes. In Sect.~\ref{sect_gl}
we treat the analog of Sect.~\ref{sect_classic} for
the general linear Lie algebras.
\section{Null cones associated to $O_n$ and $Sp_n$} \label{sect_classic}
 Denote by $GL_n, O_n$, and $Sp_{2n}$ the complex general linear group,
 orthogonal group, and the symplectic group respectively.
 Denote by $M_{n,m}$ the set of all complex $n \times m$ matrices.
\subsection{The orthogonal setup} \label{subsect_orth}
 Let $V$ be a vector space of complex dimension $n$ with
 a non-degenerate symmetric bilinear form. We
 identify $V$ with ${\C}^n$ endowed with the standard symmetric
 bilinear form
 \[
  (u, v) = \sum_{i =1}^n u_i v_i,
 \]
 where $u, v$ denote the $n$-tuples
 $(u_1, \ldots, u_n)$ and $(v_1, \ldots, v_n)$.

 For convenience, we can identify $M_{n,m}$ with the direct sum $({{\C}}^n)^m$
 of $m$ copies of ${\C}^n$. We may write a typical element of $({\C}^n)^m$
 as an $m$-tuple:
 \[
   A = (v_1, v_2, \ldots, v_m), \quad v_j \in {\C}^n.
 \]
 We define a set of quadratic polynomials
 $\tilde{\xi}_{ij} = \tilde{\xi}_{ji}, 1 \leq i, j \leq m$
 on $({\C}^n)^m$ by
 \[
  \tilde{\xi}_{ij}(v) = (v_i, v_j) = \sum_{a=1}^n v_{ai} v_{aj}.
 \]

 \begin{remark}  \rm
   The First Fundamental Theorem for $O_n$ (cf. \cite{H})
 says that the polynomials
 $\tilde{\xi}_{ij} = \tilde{\xi}_{ji}, 1 \leq i, j \leq m,$
 generate the algebra of $O_n$ invariant polynomials on $ (\C^n)^m$.
 \end{remark}

 We identify the space $S^2( {\C}^m)$ of second symmetric tensors as
 the space of symmetric $m \times m$ matrices. Define a map
 \begin{eqnarray*}
  Q : M_{n,m} & \longrightarrow & S^2 ({\C}^m)  \\
    T                    & \mapsto        & T^t T.
 \end{eqnarray*}
Here $T^t$ denotes the transpose of $T$.
By identifying $M_{n, m}$ with $( {\C}^n )^m $, we can equivalently
define $Q(T)$ as the $m \times m$ symmetric matrix whose $(i,j)$-th
entry is equal to $\tilde{\xi}_{ij}$.
 \subsection{The symplectic setup}  \label{subsect_sympl}
 Let $V$ be a vector space of complex dimension $n$ with
 a non-degenerate symplectic (i.e. anti-symmetric) bilinear form.
 It is well known that $n$ has to be an even number, say $2p$. We
 identify $V$ with ${\C}^{n}$ endowed with the following standard
 symplectic bilinear form
 \[
    ( v, v' ) = \sum_{j =1}^{n/2} (x_j y_j' - y_j x_j').
 \]
 Here $v = (x, y)$ and $v' = (x', y')$ are elements of ${\C}^{2p}$,
 expressed as $2$-tuples of elements of ${\C}^p$.

 We identify $M_{n,m}$ with the direct sum $({\C}^{n})^m$
  of $m$ copies of ${\C}^{n}$ as before.
  We may write a typical element of $({\C}^{n})^m$
  as an $m$-tuple:
  \[
   A = (v_1, v_2, \ldots, v_m), \quad v_j \in {\C}^{n}.
  \]
 We define a set of quadratic polynomials
 $\check{\xi}_{ij} = - \check{\xi}_{ji}, 1 \leq i, j \leq m$
 on $(\C^{n})^m$ by
 \[
  \tilde{\xi}_{ij}(v) = ( v_i, v_j ).
 \]

\begin{remark} \rm
   The First Fundamental Theorem for $Sp_n$ (cf. \cite{H})
 says that the polynomials
 $\check{\xi}_{ij} = - \check{\xi}_{ji}, 1 \leq i, j \leq m,$
 generate the algebra of $Sp_n$ invariant polynomials on $ ({\C}^n)^m$.
\end{remark}

 We identify the space $\Lambda^2( {\C}^m)$ of
second anti-symmetric tensors in $\C^m$ as
 the space of anti-symmetric $m \times m$ matrices. Denote by
 \[ J= \left[
         \begin{array}{rr}
              0 & I_p \\
           - I_p &  0
          \end{array} \right], \]
 where $I_p$ denotes the identity $p \times p$ matrix. Define a map
 \begin{eqnarray*}
  Q_{sp} : M_{n, m} & \longrightarrow & \Lambda^2 ({\C}^m)  \\
    T                    & \mapsto        & T^t J T.
 \end{eqnarray*}
 We also write $Q_{sp}$ as $Q$ when no ambiguity arises.
 By identifying $M_{n, m}$ with $( {\C}^n )^m $, we can equivalently
define $Q(T)$ as the $m \times m$ anti-symmetric matrix whose $(i,j)$-th
entry equal to $\check{\xi}_{ij}$.
\subsection{A resolution of singularities of $\NCQ$}
    Let $V$ be a vector space of complex dimension $n$
with a nondegenerate bilinear form, either symmetric or anti-symmetric.
Let $G(V)$ be the isometry group of the form, which is $O_n$
in symmetric case and $Sp_n$ in anti-symmetric case.
Let $\frak{g} (V)$ be its Lie algebra.

Given $A$ in $S^2 ({\C}^m)$ (resp. $\Lambda^2 ({\C}^m)$),
its inverse image $Q^{-1} (A)$ under the map $Q$
 is referred to as the fiber of $Q$ over $A$.
Of particular interest is the
fiber over $0$ (cf. \cite{H}),
which we will refer to as the {\em null cone} (since it is
obviously preserved by scalar dilations)
and denote by $\NCQ$. Equivalently, the variety $\NCQ$
is the set of $n \times m$ matrices on which all
$G(V)$ invariant polynomials on $M_{n,m}$ take value $0$.
Often we will think of $M_{n,m}$ as the space
$\Hom (\C^m, V)$ of all linear maps
from ${\C}^m$ to $V$ (or ${\C}^n$). An element
in $\NCQ$ is called a {\em null mapping}.
We easily have
\begin{eqnarray}
  Q (g T) & =& Q(T), \quad T \in M_{n,m}, \; g \in G(V), \label{eq_inv} \\
  Q(T h)  & =& h^t Q(T) h , \quad h \in GL_m.  \label{eq_glinv}
 \end{eqnarray}

\begin{remark} \rm
   It follows from Eqs.~(\ref{eq_inv}) and (\ref{eq_glinv})
 that $G(V)$ acts on $\NCQ$ and this action commutes
 with the action of $GL_m$. The space of
 regular functions on $\NCQ$, under the induced action of
 $G(V) \times GL_m$, has a beautiful decomposition
 into $G(V) \times GL_m $-modules (cf. \cite{H}).
\end{remark}

  A subspace of $V$ is called {\em isotropic} if any two vectors
in the subspace are orthogonal to each other
with respect to the corresponding bilinear form.
We observe that $T \in M_{n,m} = \Hom (\C^m ,V)$
is a null mapping if and only if
the image of $T$, denoted by $\Im T$, form an isotropic subspace of $\C^n$.
Denote by $J_r (V)$ the set of all $r$-dimensional
isotropic subspaces of $V$. The set $J_r (V)$ is nonempty if and
only if $r \leq n/2$. It is well known that $G(V)$ acts on
$J_r (V)$ transitively (cf. \cite{H}, Appendix 3).
The dimension of $J_r(V)$ $(r \leq n/2 )$ can be
calculated to be $r(2n -3r -1)/2$.

 Note that the null cone $\NCQ$ is a singular variety defined in terms
of a finite set of quadratic equations. The first goal of this paper is
to present a canonical resolution of singularities of $\NCQ$. Our
construction is reminiscent of the Springer resolution of
singularities of the nilpotent cone in a complex semisimple Lie algebra
(cf. e.g. \cite{CG}).

 Set $r =\min (m, [n/2])$ from now on, where $[n/2]$ denotes
the integer closest to and no larger than $n/2$.
We introduce the following variety
 \[
 \NCQT = \{ (T, U) \in \NCQ \times J_r (V)
   \mid \Im T \subset U \}.
 \]
We have the following projection maps:
   \[
   \begin{array}{rcl}
           &  \NCQT &    \\
    \mu_0 \swarrow  &        &  \searrow \pi_0 \\
        \NCQ      &        & J_r (V)
   \end{array}
   \]
The diagram is $O_n \times GL_m$-equivariant, where
$GL_m$ acts trivially on $\NCQT$ and $J_r (V)$.
$\mu_0$ is proper since $\NCQT \subset \NCQ \times J_r (V)$
and $J_r (V)$ is compact. It is easy to see that
$\pi_0$ and $\mu_0$ are both surjective.

Denote by $T_J$ the tautological bundle over
$J_r (V)$:
 \[
 T_J = \{ (u, U) \in V \times J_r (V) \mid u \in U \}.
 \]
 Denote by $\bdl$ the rank $r$ trivial bundle over $J_r (V)$.
 Given $U \in J_r (V)$, the fiber of $\pi_0$ over $U$
 is canonically identified with $\Hom (\C^m, U)$.
Thus we have
\begin{eqnarray*}
  \dim \NCQT & =& mr + \dim J_r (V) \\
              & =& mr + r(2n -3r -1)/2  \\
              & =& r(2m +2n -3r -1)/2.
\end{eqnarray*}
This proves the following characterization of $\NCQT$.

\begin{proposition} \label{prop_bdl}
  The variety $\NCQT$ is isomorphic to the
  tensor product $Taut \otimes \underline{\C}^m$
  of vector bundles $Taut$ and $\underline{\C}^m$ over
  the variety $J_r (V)$ such that the following
  diagram commutes:
   \[
   \begin{array}{rcl}
             \NCQT & \cong   & T_J \bigotimes \underline{\C}^m  \\
    \pi_0\searrow  &         &  \swarrow   \\
                   & J_r (V) &
   \end{array}
   \]
 In particular, $\NCQT$ is a smooth variety of dimension
 $ r(2m +2n -3r -1)/2$, where $r = \min (m ,[n/2])$.
\end{proposition}

\begin{theorem}  \label{th_sing}
 The map $\mu_0: \NCQT \longrightarrow \NCQ$ is a resolution
 of singularities.
\end{theorem}

\begin{demo}{Proof}
  It is clear that the set $\NCQ_0$ of all null mappings of maximal
 rank which is equal to $r = \min(m, [ n/2])$ is a Zariski-open
 subvariety of $\NCQ$. Given $T \in \NCQ_0$ there exists a unique
 $U \in J_r (V)$ containing $\Im T$, namely $\Im T$
 itself. Thus the map $\mu_0 : \NCQT \longrightarrow \NCQ$
 over an open set $\mu_0^{-1} (\NCQ_0 )$ is one-to-one.
 Together with the smoothness of $\NCQT$ provided by
 Proposition~\ref{prop_bdl}, we conclude the proof.
\end{demo}

\begin{remark}  \rm
   In the case $n \geq 2m$ and so $r =m$, we easily see that
 the map $Q$ maps surjectively to the space $Sym_m$ of $m \times m$
 symmetric matrices, and it is submersive at any point $T$ in $\NCQ_0$.
 It follows that
 \begin{eqnarray*}
  \dim \NCQ & =& \dim M_{n,m} - \dim Sym_m   \\
            & =& nm - (m^2 +m )/2 = \dim \NCQT.
 \end{eqnarray*}
 This of course was also implied by Theorem~\ref{th_sing}.
\end{remark}

  It is known that $J_r (V)$ is disconnected if and only
if we are in the orthogonal case and $n =2m$
(which implies $r =m$) (cf. \cite{H}, Appendix 2).
Let us assume that we are in such a case first of all, and so
$J_r (V)$ has two smooth connected components
(cf. \cite{H}). It follows by Proposition~\ref{prop_bdl}
that $\NCQT$ also has two smooth connected components,
denoted by $\NCQT^+$ and $\NCQT^-$ respectively.
The null cone $\NCQ$ also consists of  two irreducible
components $\NCQ^{\pm}$, which are the image of
the two connected components of $\NCQT$ respectively.

  Now assume we are in the symplectic case, or in the orthogonal
case but $n \neq 2m$. Then $J_r (V)$ is connected smooth
and so is $\NCQT$. Thus $\NCQ$ is irreducible as the image of
the surjective map $\mu_0$ of the irreducible variety $\NCQT$.
\subsection{Relations with closure of conjugacy classes}
   In this subsection,we always assume that $n \geq 2m$,
and in addition $m$ is even in the orthogonal case.
We need to reformulate the definition of $\NCQ$.

Let $W$ be a vector space of complex dimension $m$
with a non-degenerate bilinear form of type opposite
to the one on $V$. We identify $\Hom (W, V) = M_{nm}$.
Given any $T \in \Hom (W, V)$ the adjoint $T^*$ is defined by
\[
  (Tw, v)_V = (w, T^*v)_W, \quad w \in W, v \in V.
\]
Here we use the subscripts to indicate to which
vector space a bilinear form belong.

\begin{remark} \rm   \label{rem_ident}
  In the setup of Subsect.~\ref{subsect_orth} and \ref{subsect_sympl}
and $W = \C^m$, we easily check by definition that
if the bilinear form on $V$ is anti-symmetric
then $T^* = - J T^t$; if the bilinear form on $V$ is symmetric,
then $T^* = T^t J$.
\end{remark}

Consider the diagram
   \[
   \begin{array}{rcl}
    \Hom (W, V)   & \stackrel{\widetilde{Q}}{\longrightarrow}
                           & \frak{g} (W) \\
    \downarrow R  &        &   \\
    \frak{g} (V)  &        &
   \end{array}
   \]
where $\widetilde{Q}$ is given by $T \mapsto T^*T$ and
$R$ by $T \mapsto T T^*$. It is easy to check that the
images of $\widetilde{Q}$ and $ R$ lie in $\frak{g} (W)$ and
$ \frak{g} (V)$ respectively. It follows from definitions and
Remark~\ref{rem_ident} that
$\widetilde{Q}^{-1} (0) = {Q}^{-1} (0) \equiv \NCQ$.

Denote by $\Cg$ the conjugacy class
of the group $G(V)$ associated to the partition
$\lambda = (2^m, 1^{n -2m})$, which is the intersection
of $\frak{g} (V)$ with the conjugacy class of $gl(V)$
associated to $\lambda$. The closure $\Cbar$ is indeed the variety of
endomorphisms in $\frak{g} (V)$ of rank no greater than $m$
and whose image is an isotropic subpace of $V$.

The variety $\NCQ = \widetilde{Q}^{-1} (0)$
appears as a special case of the variety
$Z$ in \cite{KP2}. Recall that a quotient of a $G$-variety $M$
by the group $G$ is by definition the spectrum of the algebra of
$G$-invariant regular functions on $M$.
A special case of a theorem of Kraft and Procesi
relevant to our considerations can be formulated as follows.

\begin{theorem}
  The map $R$ maps $\NCQ$ surjectively to the closure $\Cbar$ of the
 conjugacy class $\Cg$ associated to the partition
 $\lambda = (2^m, 1^{n -2m})$. Furthermore $R$ can be identified
 with the quotient map by $G(W)$ from $\NCQ$ to $\Cbar$.
\end{theorem}

Define the variety
\begin{eqnarray*}
 \Ct = \{ (g, U) \in  \Cbar \times J_m (V)  | \Im g \subset U \} ,
\end{eqnarray*}
and  denote by $p_0$ the (surjective) projection from $\Ct$
to the first factor $\Cbar$. We can identify $\Ct$ with
a vector bundle over $J_m(V)$, whose fiber over $U \in J_m (V)$ is
the vector space
\[
  F_0 = \{ g \in \frak{g} (V) | \Im g \subset U \}.
\]
We remark that for $g \in F_0$ we have
$U \subset \ker g$ and thus $g^2 =0 \in gl(V)$.
Since $p_0 : \Ct \rightarrow \Cbar$ is one-to-one over the
open subset $p_0^{-1} (\Cg)$ consisting of the endomorphisms of
maximal rank $m$, we have proved that

\begin{proposition}
  The map
  $p_0 : \Ct \rightarrow \Cg$ is a resolution of singularities.
\end{proposition}

The quotient map $R$ by $G(W)$ induces a natural
quotient map by $G(W)$ from $\NCQT$ to $\Ct$.
We have the following commutative diagram:
\[
   \begin{array}{lcl}
    \NCQT  & \stackrel{\mu_0}{\longrightarrow}  & \NCQ  \\
    \downarrow/G(W)  &        & \downarrow/G(W)   \\
    \Ct   & \stackrel{p_0}{\longrightarrow}     & \Cbar
   \end{array}
\]

\begin{remark} \rm
   In the case $n = 2m$ and $G(V) = O_n$, the $O_n$-conjugacy
class $\Cg$ associated to $\lambda = (2^m)$ is {\em very even},
and so $\Ct$ splits into two $SO_n$-conjugacy class
$\Ct^{\pm}$ (cf. \cite{KP2}). Then the above diagram
can be refined as follow:
\[
   \begin{array}{rcl}
    \NCQT^{\pm}  & \stackrel{\mu_0}{\longrightarrow}  & \NCQ^{\pm}  \\
    \downarrow  &        & \downarrow   \\
    \Ct^{\pm}   & \stackrel{p_0}{\longrightarrow}     & \Cg^{\pm}
   \end{array}
\]
where $\Ct^{\pm}$ is defined by
\[
 \Ct = \{ (g, U) \in \Cbar \times J_m^{\pm} (V) | \Im g \subset U \} ,
\]
and $J_m^{\pm} (V)$ are the two connected components of $J_m (V)$.
\end{remark}
\section{Null cones associated to $GL_n$} \label{sect_gl}
 Let $V$ be a vector space of complex dimension $n$.
Define a map
$$
 \varphi : \Hom (V, \C^s) \times \Hom (\C^m, V) \rightarrow
 \Hom (\C^m, \C^s) = M_{s, m}
$$
by $\varphi (A, B) = AB
  = (\xi_{ij})_{1\leq i \leq s, 1 \leq j \leq m}$. Hence $\xi_{ij}$
defines a polynomial function on
$\Hom (V, \C^s) \times \Hom (\C^m, V) \cong V^m \oplus (V^*)^s$.
Denote by $ {\cal N} = \varphi^{ -1} (0).$
The group $GL_s \times GL(V) \times GL_m$ acts on $\cal N$ by
$$
 (g_1, g_2, g_3). (A, B) =
  (g_1 A g_2^{-1}, g_2 B g_3^{-1} ).
$$

\begin{remark} \rm
The First Fundamental Theorem for $GL(V)$ says that
the $\xi_{ij}, 1\leq i \leq s, 1 \leq j \leq m$
generate the algebra of $GL(V)$ invariant polynomials on
$V^m \oplus (V^*)^s$.
\end{remark}

From now on we restrict ourselves to the case $n \geq s +m $.
\subsection{A resolution of singularities of $\cal N$}
  We introduce a variety
\begin{eqnarray*}
 \widetilde{\cal N}
 & =&
 \{ (A, B, U_1, U_2) \in {\cal N} \times \F
   \mid \Im B\subset U_1 \subset U_2 \subset \ker A \} \\
 & \subset & {\cal N} \times \F,
\end{eqnarray*}
where $\ker A$ denotes the kernel of $A$, and
$\F$ denotes the generalized flag variety
$$
 \F =
 \{ (U_1, U_2) \mid U_1 \subset U_2 \subset V,
     \dim U_1 =m, \dim U_2 = n -s \}.
$$
The group $GL_s \times GL(V) \times GL_m$ acts on $\Nt$ by letting
$$
  (g_1, g_2, g_3). (A, B, U_1, U_2) =
  (g_1 A g_2^{-1}, g_2 B g_3^{-1}, g_2 U_1, g_2 U_2 ).
$$
We have the following projection maps:
   \[
   \begin{array}{lcl}
           &  \Nt &    \\
    \mu \swarrow  &        &  \searrow \pi \\
       {\cal N}   &        & \F
   \end{array}
   \]
The diagram above is $GL_s \times GL(V) \times GL_m$-equivariant,
where $G(V)$ acts naturally on the generalized flag variety $\F$ while $GL_s$
and $ GL_m$ act
on $\F$ trivially. It is easy to see that $\pi $ and $\mu$
are surjective thanks to the assumption $n \geq s +m$,
and that $\mu$ is a proper map due to the compactness of $\F$.
The fiber of $\pi$ over a point
$(U_1, U_2) \in \F$ can be identified with the vector space
\begin{equation} \label{eq_fiber}
 F= \Hom (V/ U_2, \C^s) \times \Hom (\C^m, U_1).
\end{equation}
In other words, $ \Nt $ can be identified with a vector bundle
${\cal K}$ over the generalized flag variety $\F$:
$$
 {\cal K} := \underline{\C^m } \bigotimes \Taut \bigoplus
 \underline{\C^s } \bigotimes {\cal Q}^*,
$$
where $\underline{\C^s },  \underline{\C^m }$ are
trivial bundles of rank $s$ and $m$ respectively,
and $\Taut$, ${\cal Q}^*$ are respectively the tautological bundle
and the dual quotient bundle:
\begin{eqnarray*}
 \Taut      & =& \{ (u, U_1, U_2) \mid v \in U_1 \}
             \subset V \times \F  \\
{\cal Q}^*  & =& \{ (v, U_1, U_2) \mid v \in (V/ U_2)^*\}
             \subset V^* \times \F.
\end{eqnarray*}

Summarizing, we have proved that

\begin{proposition}  \label{prop_same}
  We have the following isomorphism of vector bundles
 over the generalized flag variety $\F$:
   \[
   \begin{array}{lcl}
           \Nt & \cong & \quad {\cal K}    \\
    \pi\searrow  &        &  \swarrow \\
                 &\F &
   \end{array}
   \]
 In particular, $\Nt$ is a smooth variety of dimension $sn +mn -sm$.
 \end{proposition}

\begin{demo}{Proof}
  It remains to calculate that
 \begin{eqnarray*}
   \dim \Nt & =& \dim \F + \dim F  \\
            & =& \frac12 (n^2 - (s^2 + (n -s -m)^2 + m^2) +(s^2 +m^2)\\
            & =& sn +mn -sm,
 \end{eqnarray*}
 where the equality
 $\dim \F =  \frac12 (n^2 - (s^2 + (n -s -m)^2 + m^2)$
 follows readily from the description of the generalized flag variety $\F$
 in terms of the quotient of $GL(V)$ by an appropriate parabolic subgroup.
\end{demo}

\begin{theorem}  \label{th_maintwo}
  $\mu : \Nt \rightarrow {\cal N}$ is a resolution of singularities.
\end{theorem}

\begin{demo}{Proof}
  The subset ${\cal N}_0$ of $\cal N$ consisting of pairs
 of maximal rank matrices $(A, B)$ is a Zariski-open set in $\cal N$.
 A pair $ (A, B) \in {\cal N}_0$
 uniquely determines a pair $(U_1, U_2) \in \F$
 such that $(A, B, U_1, U_2) \in \Nt$, namely
 $U_1 = \Im B, U_2 = \ker A$. This shows that $\mu$ maps
 $\mu^{-1}({\cal N}_0) \subset \Nt$ bijectively to ${\cal N}_0$.
 Together with the smoothness of $\Nt$ proved in
 Proposition~\ref{prop_same}, we have concluded the proof.
\end{demo}

\begin{remark} \rm
  Indeed it is easy to see that ${\cal N}_0$
 is a single $GL(V)$-orbit since $n \geq s +m$. An easy
 calculation shows that $\varphi$
 is submersive over any point in ${\cal N}_0 $. So we have
 \begin{eqnarray*}
  \dim {\cal N}
            & =& \dim \Hom (V, \C^s) + \dim \Hom (\C^m, V) -\dim M_{s,m} \\
            & =& sn +mn -sm = \dim \Nt
 \end{eqnarray*}
 by Proposition~\ref{prop_same}. This fact, of course was
 implied by Theorem~\ref{th_maintwo}.
\end{remark}

 When $s =m$ and thus $ n \geq 2m$, the space $\cal N$ appears as a
special case of the variety $Z$ studied in \cite{KP1}.
It is shown \cite{KP1} that the quotient of $\cal N$
by the diagonal action of $GL_m$ is the
variety of $n \times n$ matrices of square $0$
and rank at most $m$, which is the closure $\Obar$
of the conjugacy class ${\cal O}_{\lambda} \subset gl(V)$ corresponding to
the partition $\lambda = (2^m, 1^{n -2m})$. More explicitly
the quotient map from $\cal N$ to $\Obar$
is given by $(A, B) \mapsto BA$.

Denote by
$$
  \Ot =  \{(g, U_1, U_2) \mid \Im g \subset U_1 \subset U_2 \subset \ker g\}
                \subset \Obar \times \Fm ,
$$
and by $p$ the natural surjective projection from $\Ot$ to $\Obar$.
By a similar argument which leads to Proposition~\ref{prop_same},
we have the following identification
\[
   \begin{array}{lcl}
     \Ot          & \cong &  Taut \bigotimes Q^*  \\
    \searrow      &       &  \swarrow \\
                  & \Fm   &
   \end{array}
\]
which implies that $\Ot$ is smooth. Since
the natural surjective projection $p$ from
$\Ot$ to (the first factor) $\Obar$
is one-to-one over the open
set $p^{-1} ({\cal O}_{\lambda})$, it
is a resolution of singularities. Thus we have established
\begin{proposition}
  The map $p : \Ot \longrightarrow \Obar$ is a resolution
 of singularities.
\end{proposition}

\begin{remark}   \rm
  This resolution $p : \Ot \longrightarrow \Obar$
 differs from the classical one in terms of cotangent bundle
 of a grassmannian (compare with Proposition~\ref{prop_resol2}).
\end{remark}

The quotient map from $\cal N$ to $\Obar$ induces
a quotient map by $GL_m$ from $\Nt$ to $\Ot$ which makes
the following digram
\begin{eqnarray} \label{eq_diag}
   \begin{array}{lcl}
     \Nt             & \stackrel{\mu}{\longrightarrow} & {\cal N} \\
    \downarrow /GL_m &                                 &  \downarrow /GL_m  \\
     \Ot             & \stackrel{p}{\longrightarrow} &\Obar
   \end{array}
\end{eqnarray}
commutative.
\subsection{A second resolution of singularities of $\cal N$}
   We introduce a variety
\begin{eqnarray*}
 \Ntt
 & =&
 \{ (A, B, U) \in {\cal N} \times \Gr \mid \Im B\subset U \subset \ker A \} \\
 & \subset & {\cal N} \times \Gr,
\end{eqnarray*}
where $\Gr$ denotes the grassmannian of $m$-dimensional
subspaces of $V$.
The group $GL_s \times GL(V) \times GL_m$ acts on $\Ntt$ by letting
$$
  (g_1, g_2, g_3). (A, B, U) =
  (g_1 A g_2^{-1}, g_2 B g_3^{-1}, g_2 U ).
$$
We have the following projection maps:
   \[
   \begin{array}{lcl}
                    & \Ntt &    \\
    \mu_1 \swarrow  &      &  \searrow \pi_1 \\
       {\cal N}     &      & \Gr
   \end{array}
   \]
The diagram above is $GL_s \times GL(V) \times GL_t$-equivariant,
where $G(V)$ acts naturally on $\Gr$ while $GL_s, GL_t$ act
trivially on $\Gr$. It is easy to see that $\pi_1 $ and $\mu_1$
are surjective and that $\mu_1$ is a proper map.
The fiber of $\pi_1$ over a point
$U \in \Gr $ can be identified with the vector space
\[
 F_1 = \Hom (V/ U, \C^s) \times \Hom (\C^m, U).
\]
In other words, $ \Ntt $ can be identified with a vector bundle
${\cal K}_1$ over the $\Gr$:
$$
 {\cal K}_1 := \underline{\C^m } \bigotimes \Taut_1 \bigoplus
 \underline{\C^s } \bigotimes {\cal Q}_1^*,
$$
where $\underline{\C^s },  \underline{\C^m }$ are
trivial bundles of rank $s$ and $t$ respectively,
and $\Taut_1$, ${\cal Q}_1^*$ are respectively the
following tautological bundle and the dual quotient bundle:
\begin{eqnarray*}
 \Taut_1    & =& \{ (u, U) \mid v \in U \}
             \subset V \times \Gr  \\
{\cal Q}_1^*  & =& \{ (v, U) \mid v \in (V/ U)^*\}
             \subset V^* \times \Gr.
\end{eqnarray*}

Summarizing, we have proved that

\begin{proposition} \label{prop_iso}
  We have the following isomorphism of vector bundles
 over the generalized flag variety $\Gr$:
   \[
   \begin{array}{lcl}
           \Ntt    & \cong & \quad {\cal K}_1   \\
    \pi_1 \searrow &       &  \swarrow \\
                   & \Gr   &
   \end{array}
   \]
 In particular, $\Ntt$ is a smooth variety of dimension $sn +mn -sm$.
 \end{proposition}

\begin{demo}{Proof}
  It remains to calculate that
 \begin{eqnarray*}
   \dim \Nt & =& \dim \Gr + \dim F_1  \\
            & =& m (n -m) + m^2 (n -m) s \\
            & =& sn +mn -sm.
 \end{eqnarray*}
\end{demo}

Using Proposition~\ref{prop_iso},
the following theorem can now be proved in the same way
as Theorem~\ref{th_maintwo}.

\begin{theorem}
  The map $\mu_1 : \Ntt \rightarrow {\cal N}$
 is a resolution of singularities.
\end{theorem}

Now we restrict ourselves to the case $s =m$ and thus $n \geq 2m$.
Recall that the Kraft-Procesi quotient map from $\cal N$ to
$\Obar$ is given by $(A, B) \mapsto BA$.
Denote by
$$
  \Ott =  \{(g, U) \mid \Im g \subset U \subset \ker g\}
                \subset \Obar \times \Gr ,
$$
and by $p_1$ the surjective
projection from $\Ott$ to the first factor $\Obar$.
We easily have the following identification
\[
   \begin{array}{lcl}
     \Ott     & \cong &  T^* \Gr  \\
    \searrow  &       &  \swarrow \\
              & \Gr   &
   \end{array}
\]
where $T^* \Gr$ denotes the cotangent bundle over the grassmannian.

Noting that the fiber of the projection $p_1$ over
$g \in {\cal O}_{\lambda}$ consists of a single point.
Since $\Ott$ is smooth, the projection $p_1$ from
$\Ott$ to $\Obar$ is a resolution of singularities.
Thus we have established the following classical result.

\begin{proposition}   \label{prop_resol2}
   The map $p_1 : \Ott \longrightarrow \Obar$
 is a resolution of singularities.
\end{proposition}

The quotient map from $\cal N$ to $\Obar$ induces
a quotient map by $GL_m$ from $\Ntt$ to $\Ott$ which makes
the following digram
\begin{eqnarray*}
 \begin{array}{lcl}
     \Ntt            &\stackrel{\mu_1}{\longrightarrow} & {\cal N} \\
    \downarrow /GL_m &                                  & \downarrow /GL_m  \\
      \Ott           &\stackrel{p_1}{\longrightarrow} & \Obar
 \end{array}
\end{eqnarray*}
commutative.
\subsection{A third resolution of singularities of $\cal N$}
  We introduce the variety
\begin{eqnarray*}
 \Ntd
 & =&
 \{ (A, B, U) \in {\cal N} \times \Grdual
   \mid \Im B \subset U \subset \ker A \} \\
 & \subset & {\cal N} \times \Grdual,
\end{eqnarray*}
together with the surjective projection $\mu_2$ from $\Ntd$ to
the first factor $\cal N$. The variety $\Ntd$ is a
vector bundle over $\Grdual$ of total dimension $sn + mn -sm$.

\begin{theorem}
  The map $\mu_2 : \Ntd \longrightarrow \cal{N}$
 is a resolution of singularities.
\end{theorem}

Proofs of all the statements concerning $\Ntd$ and $\Otd$
are similar to those for $\Ntt$ and $\Ott$ in the previous subsection
which we omit.

Now we restrict again to the case $s =m$ and $n \geq 2m$.
Denote by
$$
  \Otd =   \{(g, U) \mid \Im g \subset U \subset \ker g\}
                \subset \Obar \times \Grdualm.
$$
We again can identify $\Otd$ as the cotangent bundle $T^* \Grdualm$.
We can show that the natural surjective projection $p_2$ from $\Otd$ to
$\Obar$ is a resolution of singularities.

The quotient map from $\cal N$ to $\Obar$ induces
a quotient map by $GL_m$ from $\Ntd$ to $\Otd$ which makes
the following digram
\begin{eqnarray*}
 \begin{array}{lcl}
     \Ntd              & \stackrel{\mu_2}{\longrightarrow} & {\cal N} \\
    \downarrow /GL_m   &                & \downarrow /GL_m  \\
      \Otd             & \stackrel{p_2}{\longrightarrow}   & \Obar
 \end{array}
\end{eqnarray*}
commutative.

The relation among the three resolutions of $\cal N$ is
shown by the following diagram:
\[
   \begin{array}{lcl}
                    &   \Nt   &                   \\
     q_1 \swarrow  &         & \searrow q_2     \\
       \Ntt         &         & \quad \Ntd        \\
    \mu_1 \searrow  &         & \swarrow \mu_2    \\
                    & \cal{N} &
   \end{array}
\]
where the morphisms $q_1$ and $q_2$ are defined
by sending $(A, B, U_1, U_2)$ to
$(A, B, U_1)$ and $(A, B, U_2)$ respectively.

\vspace{.2in}

\noindent {\bf Acknowledgment} This paper was initiated
when I was in Yale University in 1998. It is a pleasure to thank
Roger Howe for stimulating discussions and encouragement.


\frenchspacing
\begin{thebibliography}{AFMO}

\bibitem[CG]{CG} N. Chriss and V. Ginzburg,
{\em Representation theory and complex geometry},
Birkh\"auser, 1997.

\bibitem[H]{H} R. Howe,  {\em Perspectives on invariant theory:
Schur duality, multiplicity-free actions and beyond},
Schur Lect. (Tel Aviv) (1992)
1--182, Israel Math. Conf. Proc. {\bf 8}.

\bibitem[KP1]{KP1} H. Kraft and C. Procesi, {\em Closures of conjugacy
classes are normal}, Invent. Math. {\bf 53} (1979) 227--247.

\bibitem[KP2]{KP2} H. Kraft and C. Procesi, {\em On the geometry of
conjugacy classes in classical groups}, Comment. Math. Helvetici,
{\bf 57} (1982) 539--602.

\bibitem[W]{W} W. Wang, {\em Lagrangian construction of
the $(gl_n, gl_m)$ duality}, preprint, math.RT/9907154, submitted.

\end{thebibliography}
\end{document}